\documentclass[runningheads]{llncs}
\usepackage{graphicx}
\usepackage{booktabs} 
\usepackage{multirow,tabularx}
\usepackage[hidelinks]{hyperref}
\usepackage{xcolor}
\usepackage{amsfonts}
\usepackage{amsmath}
\usepackage{amssymb}
\newcommand{\powerset}{\raisebox{.15\baselineskip}{\Large\ensuremath{\wp}}}

\DeclareMathOperator*{\argmin}{arg\,min}

\DeclareMathOperator{\inX}{\in\cX}

\newcommand{\bbR}{\mathbb{R}}

\newcommand{\bB}{\mathbf{B}}
\newcommand{\bc}{\mathbf{c}}

\newcommand{\bd}{\mathbf{d}}

\newcommand{\be}{\mathbf{e}}

\newcommand{\bff}{\mathbf{f}}

\newcommand{\bo}{\mathbf{o}}

\newcommand{\br}{\mathbf{r}}

\newcommand{\bt}{\mathbf{t}}

\newcommand{\bx}{\mathbf{x}}

\newcommand{\by}{\mathbf{y}}

\newcommand{\cA}{\mathcal{A}}

\newcommand{\cC}{\mathcal{C}}

\newcommand{\cO}{\mathcal{O}}

\newcommand{\cS}{\mathcal{S}}

\newcommand{\cX}{\mathcal{X}}

\newcommand{\norm}[1]{\lVert #1\rVert} 

\makeatletter
\def\Ddots{\mathinner{\mkern1mu\raise\p@
\vbox{\kern7\p@\hbox{.}}\mkern2mu
\raise4\p@\hbox{.}\mkern2mu\raise7\p@\hbox{.}\mkern1mu}}
\makeatother

\newcommand{\beginsupplement}{%
        \setcounter{table}{0}
        \renewcommand{\thetable}{S\arabic{table}}%
        \setcounter{figure}{0}
        \renewcommand{\thefigure}{S\arabic{figure}}%
     }

\usepackage[ruled, vlined]{algorithm2e}
\SetCommentSty{mycommfont}

\begin{document}
\title{Ensuring smoothly navigable approximation sets by B\'ezier curve parameterizations in evolutionary bi-objective optimization {\\ \small -- applied to brachytherapy treatment planning for prostate cancer}}
%
\titlerunning{Smoothly navigable approximation sets in bi-objective optimization}
%
\author{S.C. Maree\inst{1} \and
T. Alderliesten\inst{2} \and
P.A.N. Bosman\inst{1}}
\authorrunning{S.C. Maree et al.}
%
\institute{Life Sciences and Health research group, Centrum Wiskunde \& Informatica, Amsterdam, The Netherlands,
\email{\{maree,peter.bosman\}@cwi.nl} \\
 \and
Department of Radiation Oncology, Leiden University Medical Center, \\ Leiden, The~Netherlands,
\email{t.alderliesten@lumc.nl}}
\maketitle              
\begin{abstract}
The aim of bi-objective optimization is to obtain an approximation set of (near) Pareto optimal solutions. A decision maker then navigates this set to select a final desired solution, often using a visualization of the approximation front. The front provides a navigational ordering of solutions to traverse, but this ordering does not necessarily map to a smooth trajectory through decision space. This forces the decision maker to inspect the decision variables of each solution individually, potentially making navigation of the approximation set unintuitive. In this work, we aim to improve approximation set navigability by enforcing a form of smoothness or continuity between solutions in terms of their decision variables. Imposing smoothness as a restriction upon common domination-based multi-objective evolutionary algorithms is not straightforward. Therefore, we use the recently introduced uncrowded hypervolume (UHV) to reformulate the multi-objective optimization problem as a single-objective problem in which parameterized approximation sets are directly optimized. We study here the case of parameterizing approximation sets as smooth B\'ezier curves in decision space. We approach the resulting single-objective problem with the gene-pool optimal mixing evolutionary algorithm (GOMEA), and we call the resulting algorithm BezEA. We analyze the behavior of BezEA and compare it to optimization of the UHV with GOMEA as well as the domination-based multi-objective GOMEA. We show that high-quality approximation sets can be obtained with BezEA, sometimes even outperforming the domination- and UHV-based algorithms, while smoothness of the navigation trajectory through decision space is guaranteed.

\keywords{evolutionary algorithm \and multi-objective optimization \and hypervolume \and B\'ezier curve estimation \and approximation set navigation }
\end{abstract}
\section{Introduction}
The aim of multi-objective optimization is to obtain a set of solutions that is as close as possible to the set of Pareto-optimal solutions, with different trade-offs between the objective functions. A decision maker can then navigate the obtained set, called the \textit{approximation set}, to select a desired solution. The decision maker often incorporates external factors in the selection process that are not taken into account in the optimization objectives. An inspection of the decision variables of individual solutions is therefore required to determine their desirability. To guide the selection in bi-objective optimization, a visualization of the \textit{approximation front} (i.e., the approximation set mapped to objective space) or trade-off curve can be used. The approximation front then intuitively implies a navigational order of solutions by traversing the front from one end to the other. However, solutions with similar objective values could still have completely different decision values. The decision values of all solutions then need to be inspected individually and carefully because they may not change predictably when the approximation front is traversed. This could make navigation of the approximation set unintuitive and uninsightful.
\vspace{-0.02cm}

Population-based multi-objective evolutionary algorithms (MOEAs) have successfully been applied to real-world black-box optimization problems, for which the internal structure is unknown, or too complex to exploit efficiently by direct problem-specific design \cite{Deb01book,zitzler01,bouter17GECCOmogomea}. However, imposing a form of smoothness or continuity in terms of decision variables between solutions in the approximation set as a restriction upon the population of MOEAs is not straightforward. An underlying requirement to do so is that control over approximation sets as a whole is needed. However, typical dominance-based EAs use single-solution-based mechanics. Alternatively, multi-objective optimization problems can be formulated as a higher-dimensional single-objective optimization problem by using a quality indicator that assigns a fitness value to approximation sets. An interesting quality indicator is the hypervolume measure \cite{Zitzler99}, as it is currently the only known Pareto-compliant indicator, meaning that an approximation set of given size with optimal hypervolume is a subset of the Pareto set \cite{Knowles2002,Zitzler2003,Fleischer2003}. However, the hypervolume measure has large drawbacks when used as quality indicator in indicator-based optimization, as it does not take dominated solutions into account. The uncrowded distance has been recently introduced to overcome this \cite{Toure19}, which then resulted in the uncrowded hypervolume (UHV) measure \cite{Maree20ECJ}. The UHV can be used directly as a quality indicator for indicator-based multi-objective optimization. To be able to optimize approximation sets in this approach, fixed-size approximation sets are parameterized by concatenating the decision variables of a fixed number of solutions \cite{Wang17,Beume07,Maree20ECJ}. A single-objective optimizer can then be used to directly optimize approximation sets. The resulting single-objective optimization problem is however rather high-dimensional. To efficiently solve it, the UHV gene-pool optimal mixing evolutionary algorithm (UHV-GOMEA) \cite{Maree20ECJ}, exploits grey-box properties of the UHV problem by only updating a subset of the decision variables corresponding to one (or a few) multi-objective solutions.
\vspace{-0.02cm}

In this work, we go beyond an unrestricted concatenation of the decision variables of solutions and we propose to model approximation sets as sets of points that lie on a B\'ezier curve \cite{Gallier99} in decision space. Optimizing only the control points of the B\'ezier curve, that define its curvature, enforces the decision variables of solutions in the approximation set to vary in a smooth, continuous fashion, thereby likely improving intuitive navigability of the approximation set. Previous work on parameterizations of the approximation set has been applied mainly in a post-processing step after optimization, or was performed in the objective space \cite{Kobayashi18,Bhardwaj14,Mehta14}, but this does not aid in the navigability of the approximation set in decision space. Moreover, fitting a smooth curve through an already optimized set of solutions might result in a bad fit, resulting in a lower-quality approximation set. Additionally, we will show that specifying solutions as points on a B\'ezier curve directly enforces a form of diversity within the approximation set, which can actually aid in the optimization process, and furthermore reduces the problem dimensionality of the single-objective problem.

The remainder of this paper is organized as follows. In Section~\ref{sec:ppsn20_moo}, we introduce preliminaries on UHV-based multi-objective optimization.  In Section~\ref{sec:ppsn20_smoothness_measure}, we define a measure for navigational smoothness of approximation sets. In Section~\ref{sec:ppsn20_bezier}, we introduce B\'ezier curves and the corresponding optimization problem formulation. Empirical benchmarking on a set of benchmark problems is performed in Section~\ref{sec:ppsn20_experiments}. 
Finally, we discuss the results and conclude in Section~\ref{sec:ppsn20_discussion}. 
Additionally, in the Supplementary (Section~\ref{sec:bbez20_brachy}), we demonstrate BezEA on a real-world optimization problem that arises in the treatment of prostate cancer and analyze the resulting approximation sets.

\section{UHV-based multi-objective optimization}
\label{sec:ppsn20_moo}
Let $\bff : \cX \rightarrow \bbR^m$ be a to-be-minimized $m$-dimensional vector function and $\cX \subseteq \bbR^n$ be the $n$-dimensional (box-constrained) decision space. When the objectives in $\bff$ are conflicting, no single optimal solution exists, but the optimum of $\bff$ can be defined in terms of \textit{Pareto optimality} \cite{Knowles2006}. A solution $\bx\inX$ is said to \textit{weakly dominate} another solution $\by\inX$, written as $\bx\preceq\by$, if and only if $f_i(\bx) \leq f_i(\by)$ for all $i$. When the latter relation is furthermore strict (i.e., $f_i(\bx) < f_i(\by)$) for at least one $i$, we say that $\bx$ \textit{dominates} $\by$, written as $\bx\prec\by$. 
A solution that is not dominated by any other solution in $\cX$ is called \textit{Pareto optimal}. The \textit{Pareto set} $\cA^\star$ is the set of all Pareto optimal solutions, i.e., $\cA^\star = \{ \bx\inX : \nexists \by\inX : \by\prec\bx \} \subset \cX$. The image of the Pareto set under $\bff$ is called the \textit{Pareto front}, i.e., $\{ \bff(\bx) : \bx\in\cA^\star \} \subset \bbR^m$. The aim of multi-objective optimization is to approximate the Pareto set with a set of non-dominated solutions called an \textit{approximation set} $\cA$. Let $\cS \subseteq \cX$ be a solution set, that can contain dominated solutions and let $A : \powerset(\cX) \rightarrow \powerset(\cX)$ be the approximation set given by $\cS$, i.e., $A(\cS) = \{\bx\in\cS : \nexists \by\in\cS : \by\prec\bx\}$, where $\powerset(\cX)$ is the powerset of $\cX$.

\begin{figure}
\begin{center}
\includegraphics[width=\textwidth]{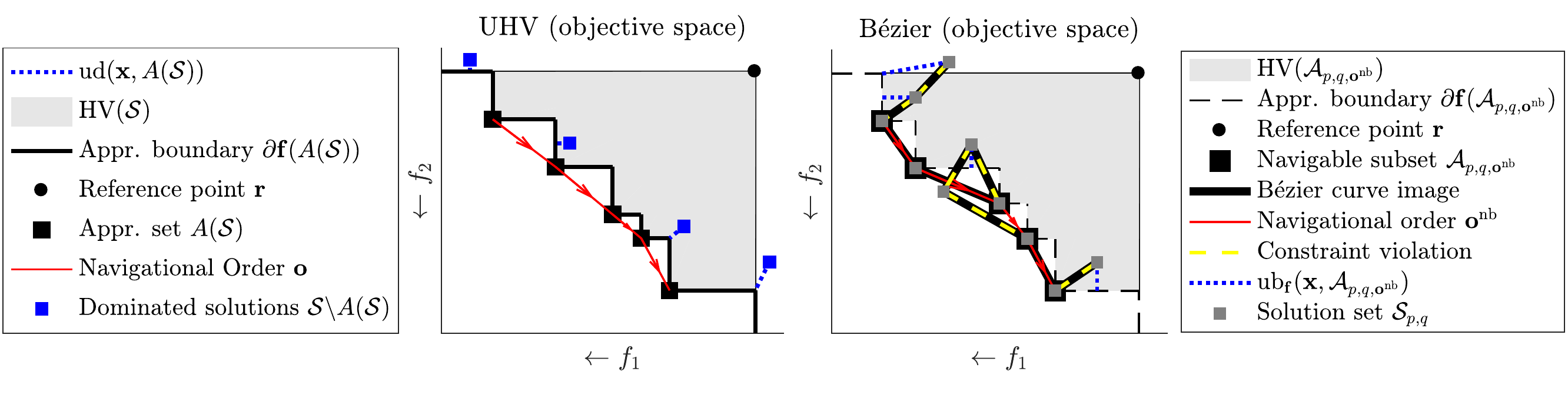}
\vspace*{-1.0cm}
\caption{Illustration of the uncrowded hypervolume (UHV) \cite{Maree20ECJ} (left) for a bi-objective minimization problem, and the B\'ezier parameterization (right).}
\vspace*{-0.2cm}
\label{fig:ppsn20_uhv}
\end{center}
\end{figure}

The hypervolume measure $\mbox{HV} : \powerset(\cX) \rightarrow \bbR$ \cite{Zitzler2003,auger05} measures the area or volume dominated by all solutions in the approximation set, bounded by a user-defined reference point $\mathbf{r}\in\bbR^m$, as shown in Figure~\ref{fig:ppsn20_uhv}. As the hypervolume ignores dominated solutions, we use the \textit{uncrowded distance} to assign a quality value to dominated solutions \cite{Toure19}. The uncrowded distance $\text{ud}_\bff(\bx,\cA)$ measures the shortest Euclidean distance between $\bx$ and the \textit{approximation boundary} $\partial \bff(\cA)$, when $\bx$ is dominated by any solution in $\cA$ or outside the region defined by $\br$, and is defined $\text{ud}_\bff(\bx,\cA) = 0$ else (Figure~\ref{fig:ppsn20_uhv}). It is called the uncrowded distance as the shortest distance to $\partial \bff(\cA)$ is obtained for a point on the boundary that is not in $\cA$ itself. Combining the uncrowded distance with the hypervolume measure results in the \textit{uncrowded hypervolume} (UHV) \cite{Maree20ECJ},
\vspace*{-0.2cm}
\begin{equation}
\label{eqn:ppsn20_uhv}
\text{UHV}_\bff(\cS) = \text{HV}_\bff(\cS) - \frac{1}{|\cS|} \sum_{\bx\in\cS} \text{ud}_\bff(\bx, A(\cS))^m.
 \vspace*{-0.2cm}
\end{equation}
We use the subscript $\bff$ to denote that its value is computed with respect to the multi-objective problem $\bff$.  To be able to optimize the UHV of a solution set, a parameterization of solution sets is required. Let $\phi\in\bbR^l$ be such a parameterization consisting of $l$ decision variables, and let $S(\phi) = \{\bx_1,\bx_2,\ldots\}$ be an operator that transforms $\phi$ into its corresponding solution set. The resulting UHV-based optimization problem is then given by, 
\vspace*{-0.2cm}
\begin{equation}
\setlength{\jot}{0pt} 
\label{eqn:ppsn20_UHVmop}
\begin{split}
\text{maximize} \quad & \text{UHV}_{\bff,S}(\phi) 
= \text{HV}_{\bff}(S(\phi)) -  \frac{1}{|\cS(\phi)|} \sum_{\bx\in\cS(\phi)} \text{ud}_\bff(\bx, A(\cS(\phi)))^m, 
\\
\text{with} \quad & \bff : \cX \subseteq \bbR^n \rightarrow \bbR^m, \quad S : \mathbb{R}^l \rightarrow \powerset(\cX), \quad \phi \in \mathbb{R}^l.
\end{split}
\vspace*{-0.2cm}
\end{equation}
In a parameterization that is commonly used, solution sets $\cS_p$ of fixed size $p$ are considered, and the decision variables of the solutions in $\cS_p$ are simply concatenated, i.e., $\phi = [\bx_1 \cdots \bx_p] \in \bbR^{p\cdot n}$ \cite{Wang17,Beume07,Maree20ECJ}. Using this parameterization, the resulting single-objective optimization problem is $l = p\cdot n$ dimensional. In \cite{Maree20ECJ}, GOMEA \cite{bouter17GECCOgomea} was used to efficiently solve this problem by exploiting the \textit{grey-box} (gb) property that not all solutions $\bx_i$ have to be recomputed when only some decision variables change. The resulting algorithm, which we call UHVEA-gb here (and was called UHV-GOMEA-Lm in \cite{Maree20ECJ}), greatly outperformed the mostly similar algorithm UHVEA-bb (called UHV-GOMEA-Lf in \cite{Maree20ECJ}) but in which the UHV was considered to be a \textit{black box} (bb). This problem parameterization however does not guarantee any degree of navigational smoothness of the approximation set, which is the key goal in this paper. 


\section{A measure for navigational smoothness}
\label{sec:ppsn20_smoothness_measure}
We introduce a measure for the navigational smoothness of an approximation set. Let $\cS_p = \{\bx_1,\bx_2,\ldots,\bx_p\}$ be an approximation set of size $p$. Furthermore, let  the \textit{navigation order} $\bo$ be a permutation of (a subset of) $I = \{1,2,\ldots,p\}$, representing the indices of the solutions in $\cS_p$ that the decision maker assesses in the order the solutions are inspected. The \textit{(navigational) smoothness} $\mbox{Sm}(\cS_p,\bo)$ is then defined as,
\vspace*{-0.2cm}
\begin{equation}
\mbox{Sm}(\cS_p,\bo) = \frac{1}{p-2}\sum_{i = 2}^{p-1} \frac{ \norm{\bx_{\bo_{i-1}} - \bx_{\bo_{i+1}}}}{\norm{\bx_{\bo_{i-1}} - \bx_{\bo_i}} + \norm{\bx_{\bo_{i}} - \bx_{\bo_{i+1}}} }.
 \vspace*{-0.2cm}
\end{equation} 
This smoothness measure measures the \textit{detour length}, i.e., the extra distance traveled (in decision space) when going to another solution via an intermediate solution, compared to directly going there. 

Throughout this work, we will consider a navigational order $\bo$ for approximation sets $\cA$ such that $f_1(\bx_{\bo_i}) < f_1(\bx_{\bo_j})$ holds whenever $i < j$ holds, i.e., from left to right in the objective space plot Figure~\ref{fig:ppsn20_uhv}. We therefore simply write $\mbox{Sm}(\cA,\bo) = \mbox{Sm}(\cA)$ from now on. Note that $\mbox{Sm}(\cA) \in [0,1]$, and only if all solutions are colinear in decision space, $\mbox{Sm}(\cA) = 1$ holds. This we consider the ideal scenario, where the decision variables of solutions change perfectly predictably. This also implies that any other (continuous) non-linear curve is not considered to be perfectly smooth. Although one could argue for different definitions of smoothness, we will see later that this measure serves our purpose for distinguishing smoothly from non-smoothly navigable approximation sets.

\section{B\'ezier curve parameterizations of approximation sets}
\label{sec:ppsn20_bezier}
\begin{figure}[t]
\begin{center}
\includegraphics[width=\textwidth]{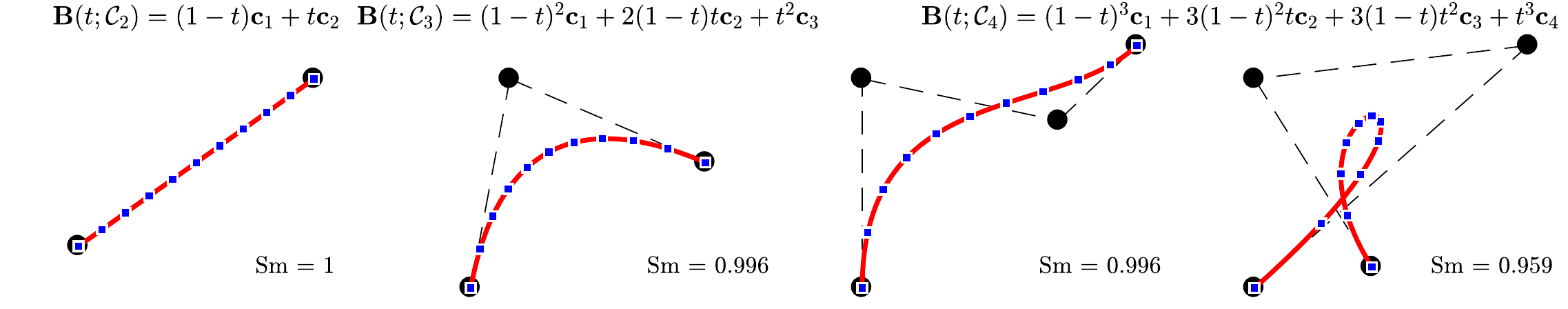}
\end{center}
\vspace*{-0.6cm}
\caption{Illustration of B\'ezier curves (red) in decision space with different control points (black). Blue points correspond to $p = 10$ evenly spread values of $t$, and the smoothness (Sm) of these $p$ points  is given, computed based on $\bo^\text{bez}$.}
\vspace*{-0.3cm}
\label{fig:ppsn20_bezier}
\end{figure}

A B\'ezier curve $\bB(t ; \cC_q)$ is a parametric curve that is commonly used in computer graphics and animations to model smooth curves and trajectories \cite{Gallier99}. An $n$-dimensional B\'ezier curve is fully specified by an ordered set of $q \geq 2$ control points $\cC_q = \{\bc_1,\ldots,\bc_q\}$ with  $\bc_j\in\bbR^n$, and given by,
\begin{equation}
\label{eqn:ppsn20_bezier}
\bB(t ; \cC_q)= \sum_{j = 1}^q b_{j-1,q-1}(t) \bc_j, \quad \mbox{with} \quad b_{j,q}(t) := {q \choose j}(1-t)^{q-j}t^{j},
\end{equation}
for $0\leq t\leq1$, where ${q \choose j}$ are the binomial coefficients. Examples of B\'ezier curves are shown in Figure~\ref{fig:ppsn20_bezier}. The first and last control points are always the end points of the B\'ezier curve, while intermediate control points do not generally lie on the curve. 
We parameterize a solution set $\cS_p = \{\bx_1,\ldots, \bx_p\}$ of fixed size $p$ using an $n$-dimensional B\'ezier curve $\bB(t ; \cC_q)$ with $q$ control points. On this curve, $p$ points $\bx_{i} = \bB\left({( i - 1 )}/{( p - 1 ) }; \cC_q\right)$ are selected, evenly spread in the domain of $t$. 
The resulting solution set $S_{p,q}(\phi) = \{ \bx_1, \bx_2,\ldots,\bx_p\}$ is then given by,
$$S_{p,q}(\phi) = \left\{ \bB\left(\frac{0}{p - 1 }; \cC_q\right),  \bB\left( \frac1{ p - 1}; \cC_q\right), \ldots, \bB\left(\frac{p-1} { p - 1 }; \cC_q\right) \right\},$$
with $\phi = [\bc_1 \cdots \bc_q] \in \bbR^{q\cdot n}$. Note that inverting the order of control points does not affect the B\'ezier curve. To avoid this symmetry in the parameterization, we standardize the curve direction throughout optimization. After a change of the curve, we check if $f_1(\bc_1) < f_1(\bc_q)$ holds. If not, the order of the control points is simply inverted. 

\vspace{-0.5cm}
\begin{algorithm}
\SetKwInOut{Function}{function}
\SetKwInOut{Input}{input}
\SetKwInOut{Output}{output}
\Function{[$\cA_{p,q,\bo^\text{nb}}, ( \bo^\text{nb})] = \mbox{A}^\text{nb}(\cS_{p,q}, \bo^\text{bez})$}
\Input{B\'ezier solution set $\cS_{p,q} = \{\bx_1,\ldots,\bx_p\}$ with intrinsic ordering $\bo^\text{bez}$}
\Output{Approximation (sub)set $\cA_{p,q,\bo^\text{nb}}$, (navigational order $\bo^\text{nb}$), }
\BlankLine
$\eta = \argmin_{i\in\{1,\ldots,p\}}f_1(\bx_{\bo^\text{bez}_i})$\;
$\bo^\text{nb} = [ \bo^\text{bez}_\eta ]$ and $\cA_{p,q,\bo^\text{nb}} = \{ \bx_{\bo^\text{bez}_\eta} \}$\;
\For{$j = \eta,\ldots,p$ }
{
\If{$\bx_{\bo^\text{\normalfont bez}_j} \in A(\cS_{p,q})$ \normalfont{\textbf{and}} $f_2(\bx_{\bo^\text{\normalfont bez}_j}) < f_2(\bx_{\bo^\text{\normalfont nb}_\text{\normalfont end}})$}{
$\bo^\text{nb} = [ \bo^\text{nb} \;;\; \bo^\text{bez}_j]$ and $\cA_{p,q,\bo^\text{nb}} = \cA_{p,q,\bo^\text{nb}} \cup \{ \bx_{\bo^\text{bez}_j} \}$\tcp*{here $\bo^\text{nb}_\text{end} = \bo^\text{bez}_j$}
}
}
\caption{Navigational order for B\'ezier parameterizations}
\label{alg:onb}
\end{algorithm}
\vspace*{-1cm}

\subsection{A navigational order for B\'ezier parameterizations}
Solution sets $\cS_{p,q} = S_{p,q}(\phi)$ parameterized by a B\'ezier curve introduce an intrinsic order $\bo^\text{bez}$ of solutions by following the curve from $t = 0$ to $t = 1$. Even though the solutions in $\cS_{p,q}$ now lie on a smooth curve in decision space, it might very well be that some of these solutions dominate others. We define a navigational-B\'ezier (nb) order $\bo^\text{nb}$ for a solution set $\cS_{p,q}$ that follows the order of solutions $\bo^\text{bez}$ along the B\'ezier curve, but also aligns with the left-to-right ordering described in Section~\ref{sec:ppsn20_smoothness_measure}. Pseudo code for $\bo^\text{nb}$ is given in Algorithm~\ref{alg:onb}, and an example is given in Figure~\ref{fig:ppsn20_uhv}. The navigational order $\bo^\text{nb}$ starts from the solution with best $f_1$-value and continues to follow the B\'ezier curve (i.e., in the order $\bo^\text{bez}$) until the solution with best $f_2$-value is reached, only improving in $f_2$ (and thereby worsening in $f_1$) along the way, and skipping solutions that violate this property. Let $\cA_{p,q,\bo^\text{nb}} = A^\text{nb}(\cS_{p,q},\bo^\text{bez})$ be the resulting subset of $\cS_{p,q}$ pertaining to exactly the solution indices as specified in $\bo^\text{nb}$, and note that this is an approximation set.

\subsection{Unfolding the B\'ezier curve (in objective space)}
Smoothly navigable approximation sets can now be obtained by maximizing the hypervolume of $\cA_{p,q,\bo^\text{nb}}$. To maximize the number of navigable solutions $|\cA_{p,q,\bo^\text{nb}}| = |\bo^\text{nb}|$, we need to unfold the B\'ezier curve in objective space. For this, we introduce a constraint violation function $C(\cS_{p,q},\bo^\text{nb}) \geq 0$, as given in Algorithm~\ref{alg:bezier_constraint} and illustrated in Figure~\ref{fig:ppsn20_uhv}. It is composed of two parts. The first part is similar to the uncrowded distance term in Eqn.~\eqref{eqn:ppsn20_uhv}, but the approximation boundary is now given by $\cA_{p,q,\bo^\text{nb}}$. The second part aims to pull solutions that are not in $\cS_{p,q,\bo^\text{nb}}$ towards neighboring solutions on the B\'ezier curve.

\vspace*{-0.5cm}
\begin{algorithm}
\SetKwInOut{Function}{function}
\SetKwInOut{Input}{input}
\SetKwInOut{Output}{output}
\Function{$C(\cS_{p,q},\bo^\text{bez}) \geq 0$}
\Input{B\'ezier solution set $\cS_{p,q} = \{\bx_1,\ldots,\bx_p\}$ with intrinsic ordering $\bo^\text{bez}$}
\Output{Constraint value $C \geq 0$}
\BlankLine
[$\cA, \bo^\text{nb}] =  \mbox{A}^\text{nb}(\cS_{p,q}, \bo^\text{bez}$)\tcp*{See Algorithm~\ref{alg:onb}}
$ C = \frac{1}{|\cS_{p,q}|} \sum_{\bx\in\cS_{p,q}} \text{ud}_\bff(\bx, \cA^m)$\tcp*{Uncrowded distance (ud), see \eqref{eqn:ppsn20_uhv}}
\BlankLine
\For{$j = 1,\ldots,|\cS_{p,q}|-1$}
{
\If{$\bo^\text{\normalfont bez}_{j} \notin \bo^\text{nb}  $ or $\bo^\text{\normalfont bez}_{j+1} \notin \bo^\text{nb}$}
{
$C = C + \norm{\bff(\bx_{\bo^\text{bez}_{j}}) -\bff(\bx_{\bo^\text{bez}_{j+1}}) }$\tcp*{Euclidean distance in $\bbR^m$}
} 
} 
\caption{B\'ezier constraint violation function}
\label{alg:bezier_constraint}
\end{algorithm}
\vspace*{-1cm}

\subsection{B\'ezier parameterization + GOMEA = BezEA}
\label{sec:ppsn20_gomea}
The resulting B\'ezier curve optimization problem is given by,
\begin{equation}
\vspace{-0.2cm}
\label{eqn:ppsn20_bezierMOP}
\setlength{\jot}{0pt} 
\begin{split}
\text{maximize} \quad & \text{HV}_{\bff,S_{p,q}}(\phi) = \text{HV}_\bff(A^\text{nb}(S_{p,q}(\phi))), 
\\
\text{with} \quad & C(S_{p,q}(\phi),\bo^\text{nb}(\phi)) = 0, \\
& \bff : \cX \subseteq \bbR^n \rightarrow \bbR^m, \quad S_{p,q} : \mathbb{R}^{q\cdot n} \rightarrow \powerset(\cX), \quad \phi \in \mathbb{R}^{q\cdot n}.
\end{split}
\vspace{-0.2cm}
\end{equation}
We use \textit{constraint domination} to handle constraint violations \cite{Deb00}. With constraint domination, the fitness of a solution is computed regardless of its feasibility. When comparing two solutions, if both are infeasible (i.e., $C > 0$), the solution with the smallest amount of constraint violation is preferred. If only one solution is infeasible, the solution that is feasible is preferred. Finally, if both solutions are feasible (i.e., $C = 0$), the original ranking based on fitness is used.

B\'ezier curves have no local control property, meaning that a change of a control point affects all solutions on the curve. Partial evaluations can therefore no longer be exploited with this parameterization, and we thus solve this problem with the black-box version of GOMEA. Analogous to the UHV naming, we brand the resulting algorithm B\'ezier-GOMEA-bb, which we abbreviate to BezEA. A detailed description of GOMEA can be found in \cite{bouter17GECCOgomea}, and a description of UHV-GOMEA in \cite{Maree20ECJ}.

\section{Numerical Experiments}
\label{sec:ppsn20_experiments}
We compare BezEA with UHVEA-gb and UHVEA-bb. These methods use a different hypervolume-based representation of the multi-objective problem, but use very similar variation and selection mechanisms, making the comparison between these methods most fair. We use the guideline setting for the population size $N$ of GOMEA with full linkage models in a black-box setting \cite{bosman13}, which for separable problems yields $N = \lfloor 10\sqrt{l}\rfloor$ and for non-separable problems $N = 17 + \lfloor 3l^{1.5}\rfloor$. BezEA solves a single-objective problem of $l = qn$ decision variables. UHVEA-bb solves a single objective problem of $l = pn$ decision variables. UHVEA-gb solves the same problem by not considering all $pn$ decision variables simultaneously, but by updating only subsets of $l = n$ decision variables, on which we base the population size guideline for UHVEA-gb. 

We furthermore include the domination-based MO-GOMEA \cite{bouter17GECCOmogomea}. In MO-GOMEA, a population of solutions is aimed to approximate the Pareto front by implicitly balancing diversity and proximity. From a population of $N_{mo}$ solutions, truncation selection is performed based on domination rank. The resulting selection is clustered into $K_{mo}$ overlapping clusters that model different parts of the approximation front. For each cluster, a Gaussian distribution is estimated to sample new solutions from, which uses very similar update rules as the single-objective GOMEA, and therefore allows for a most fair comparison to BezEA and UHVEA. MO-GOMEA obtains an elitist archive, aimed to contain 1000 solutions. For a fair comparison to the hypervolume-based methods that obtain an approximation set of at most $p$ solutions, we reduce the obtained elitist archive of MO-GOMEA to $p$ solutions using greedy hypervolume subset selection (gHSS) \cite{Guerreiro16}, which we denote by MO-GOMEA*. As described in \cite{Maree20ECJ}, to align MO-GOMEA with the other algorithms, we set $N_{mo} = p\cdot N$ and $K_{mo} = 2p$ such that the overall number of solutions in the populations is the same, and all sample distributions are estimated from the same number of solutions. 

As performance measure, we define $\Delta\text{HV}_p = \text{HV}_p^\star - \text{HV}(\cA_p)$ as the distance to the optimal hypervolume $\text{HV}_{p}^\star$ obtainable with $p$ solutions, empirically determined with UHVEA. 

\subsection{Increasing $q$}
We illustrate how increasing the number of control points $q$ of the B\'ezier curve improves achievable accuracy of BezEA (with $q = \{2,\ldots,10\}$ and $p = 10$) in case the Pareto set is non-linear. For this, we construct a simple two-dimensional problem \textit{curvePS}, with objective functions $f_1^\text{curvePS}(\bx) = (x_1 -1)^2 + 0.01 x_2^2 $ and $f_2^\text{curvePS}(\bx) = x_1^2 + (x_2-1)^2$. A large computational budget was used to show maximally achievable hypervolume, and standard deviations are therefore too small to be visible.

\begin{figure}
\begin{center}
\includegraphics[width=\textwidth]{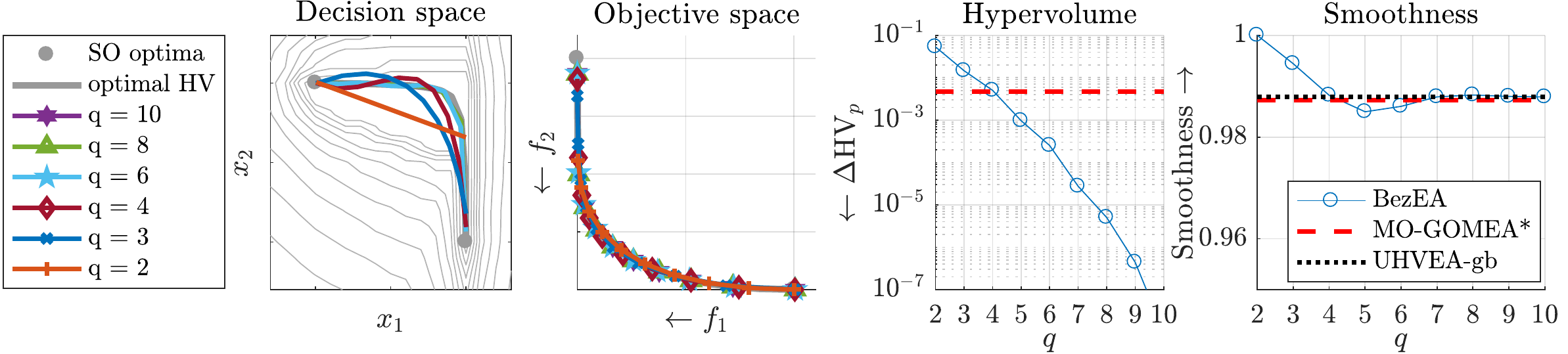}
\end{center}
\vspace*{-0.6cm}
\caption{ B\'ezier curve approximations of the Pareto set of the curvePS problem (left), obtained with BezEA. Contour lines show domination ranks, the corresponding approximation fronts (middle), and $\Delta\mbox{HV}_{10}$ together with smoothness (right).}
\label{fig:ppsn20_q}
\end{figure}

Results are shown in Figure~\ref{fig:ppsn20_q}. 
A larger $q$ results in a better approximation of the leftmost endpoint of the Pareto front (second subfigure), thereby improving $\Delta\mbox{HV}_p$ (third subfigure), but slightly lowering smoothness (fourth subfigure), as the B\'ezier curve deviates from a straight line. MO-GOMEA*, UHVEA-gb, and BezEA for large $q$ all obtain a very similar smoothness. As MO-GOMEA* does not explicitly optimize the hypervolume of its approximation set, it obtains a slightly different distribution of solutions, which results in a lower hypervolume. Additionally, MO-GOMEA* does not converge to the Pareto set due to the finite population size and inifitely large Pareto set, as described in more detail in \cite{Maree20ECJ}. Even though this is a fundamental limitation of domination-based MOEAs, this level of accuracy is often acceptable in practice.

\subsection{Comparison with UHV optimization}
Next, we demonstrate the behavior of BezEA compared to UHVEA on the simple \textit{bi-sphere} problem, which is composed of two single-objective sphere problems, $f_\text{sphere}(\bx) = \sum_{i = 1}^n x_i^2,$ of which one is translated, $f^\text{bi-sphere}_1(\bx) = f_\text{sphere}(\bx)$,  and $f^\text{bi-sphere}_2 = f_\text{sphere}(\bx - \be_1),$ where $\be_i$ is the $i^\text{th}$ unit vector. 
We set $n = 10$, and initialize all algorithms in $[-5,5]^n$. This is a separable problem and we therefore use the univariate population size guideline (i.e., $N = \sqrt{l}$). We consider the cases $p = \{10,100\}$. The computational budget is set to $2p\cdot 10^4$ evaluations of the multi-objective problem given by $\bff$ (MO-fevals). When the desired number of solutions $p$ along the front is large, neighboring solutions are nearby each other on the approximation front. This introduces a dependency between these solutions, which needs to be taken into account in the optimization process to be able to effectively solve the problem \cite{Maree20ECJ}.

\begin{figure}[t]
\includegraphics[width=\textwidth]{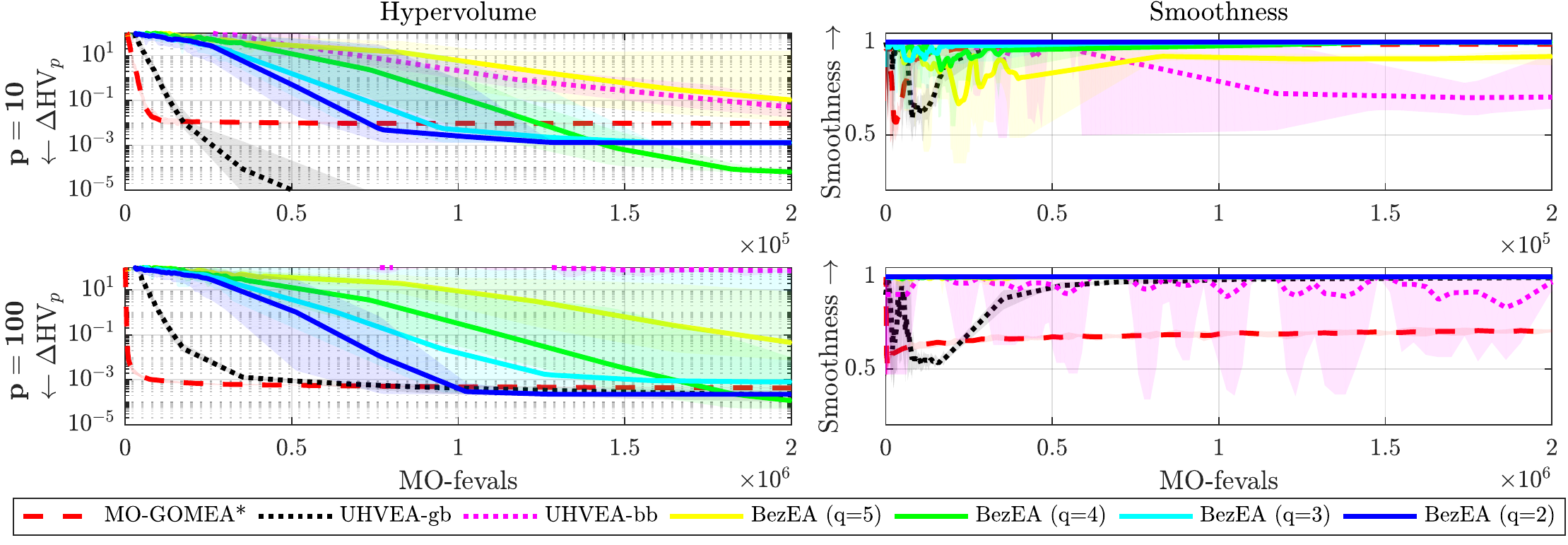}
\vspace*{-0.7cm}
\caption{Comparison of UHVEA with BezEA and MO-GOMEA* on the bi-sphere problem with $n = 10$ and $p = 10$ (top row) and $p = 100$ (bottom row). Left two subfigures show mean scores, and the shaded areas represent min/max scores, obtained over 10 runs. Objective and decision space subfigures show results of a single run. Solutions in the decision space projection are sorted based on their $f_0$-value, from best to worst.}
\label{fig:ppsn20_overhead}
\vspace*{-0.3cm}
\end{figure}

Results are shown in Figure~\ref{fig:ppsn20_overhead}. This problem is unimodal with a linear Pareto set, and the smoothness of (a subset) of the Pareto set is therefore 1.0. As UHVEA-gb converges to a subset of the Pareto set (see \cite{Maree20ECJ}), it ultimately obtains a smoothness of 1.0, even though its smoothness is initially lower. MO-GOMEA* does not converge to the Pareto set, and its smoothness stagnates close to 1.0 when $p = 10$, but stagnates around 0.7 when $p = 100$. BezEA with $q = 2$ has per construction a perfect smoothness of 1.0, and for $q = 3$ and $q = 4$, the obtained smoothness is close to 1. With $q = 5$ control points, BezEA does not converge within the given budget, resulting in a lower smoothness within the computational budget. UHVEA-gb furthermore shows a better convergence rate, which could be because UHVEA-gb can exploit partial evaluations, while this is not possible with BezEA. However, UHVEA-bb, which also does not perform partial evaluations, is unable to solve the problem for $p = 100$. This difference between BezEA and UHVEA-bb could be attributed to the lower degree of freedom that BezEA has due to the rather fixed distribution of solutions. This distribution does however not exactly correspond to the distribution of $\mbox{HV}_p^\star$. This is why a stagnation in terms of hypervolume convergence can be observed for small values of $q$. The solutions of BezEA are equidistantly distributed along the curve in terms of $t$. By doing so, intermediate control points can be used to adapt the distribution of solutions (when $q > 2$). This is why BezEA with $q = 4$ can obtain a better $\Delta\mbox{HV}_p$ than BezEA with $q = 2$, even though the Pareto set is linear. For $p = 100$, BezEA obtains a better $\Delta\mbox{HV}_p$ than UHVEA-gb, which can be explained by the increased problem complexity when the desired number of solutions along the front is large. Increasing the population size $N$ of UHVEA-gb would (at least partially) overcome this, but we aimed here to show that BezEA does not suffer from this increased complexity as its problem dimensionality depends on $q$, not $p$.

\subsection{WFG benchmark}
We benchmark BezEA, UHVEA, and MO-GOMEA on the nine commonly used WFG functions \cite{Huband2005}. We consider bi-objective WFG problems with $n = 24$ decision variables of which $k_\text{WFG} = 4$ are WFG-position variables. We furthermore set $p = 9$ and a computational budget of $10^7$ MO-fevals. A population size of $N = 200$ was shown to work well for UHVEA \cite{Maree20ECJ}, which we use here also for BezEA. We perform 30 runs, and a pair-wise Wilcoxon rank-sum test with $\alpha = 0.05$ is used to test whether differences with the best obtained result are statistically significant (up to 4 decimals). Ranks (in brackets) are computed based on the mean hypervolume values.

\begin{table}[t]
\scriptsize
\caption{Obtained hypervolume $\text{HV}_p$ (mean $\pm$ standard deviation (rank)) and mean navigational smoothness (Sm) for the 9 WFG problems with $p = 9$ solutions. Bold are best scores per problems, or those not statistically different from it.}
\label{tab:ppsn20_wfg}
\vspace*{-0.1cm}
\begin{tabular*}{\textwidth}{l@{\extracolsep{\fill}}rrrrrrrr@{\extracolsep{\fill}}}
\toprule
\# & \multicolumn{2}{r}{MO-GOMEA*} & \multicolumn{2}{r}{UHVEA-gb} & \multicolumn{2}{r}{BezEA ($q = 2$)} & \multicolumn{2}{r}{BezEA ($q = 3$)} \\ 
& $\mbox{HV}_9$ & Sm & $\mbox{HV}_9$ & Sm & $\mbox{HV}_9$ & Sm & $\mbox{HV}_9$ & Sm \\
 \midrule 
1  & $\textbf{97.60}\pm 0.7$ (1) &\textit{0.76} & $93.62 \pm 1.7$ (2) &\textit{0.67} & $90.35 \pm 1.1$ (4) &\textit{1.00} & $90.37 \pm 1.2$ (3) &\textit{0.99}\\ 
2  & $110.09 \pm 0.0$ (2) &\textit{0.86} & $\textbf{110.38}\pm 1.0$ (1) &\textit{0.66} & $97.74 \pm 0.0$ (4) &\textit{1.00} & $97.85 \pm 0.0$ (3) &\textit{0.98}\\ 
3  & $116.11 \pm 0.1$ (4) &\textit{0.93} & $116.42 \pm 0.1$ (3) &\textit{0.71} & $\textbf{116.50}\pm 0.0$ (1) &\textit{1.00} & $\textbf{116.50}\pm 0.0$ (2) &\textit{1.00}\\ 
4  & $111.88 \pm 0.8$ (3) &\textit{0.75} & $\textbf{112.37}\pm 0.7$ (1) &\textit{0.69} & $111.59 \pm 1.3$ (4) &\textit{1.00} & $\textbf{112.19}\pm 1.3$ (2) &\textit{0.98}\\ 
5  & $112.03 \pm 0.1$ (3) &\textit{0.66} & $111.86 \pm 0.3$ (4) &\textit{0.63} & $112.17 \pm 0.0$ (2) &\textit{1.00} & $\textbf{112.19}\pm 0.0$ (1) &\textit{1.00}\\ 
6  & $113.86 \pm 0.3$ (3) &\textit{0.88} & $114.23 \pm 0.2$ (2) &\textit{0.72} & $\textbf{114.34}\pm 0.1$ (1) &\textit{1.00} & $113.02 \pm 0.3$ (4) &\textit{0.99}\\ 
7  & $114.06 \pm 0.1$ (4) &\textit{0.94} & $114.32 \pm 0.1$ (3) &\textit{0.66} & $114.37 \pm 0.0$ (2) &\textit{1.00} & $\textbf{114.38}\pm 0.0$ (1) &\textit{1.00}\\ 
8  & $110.70 \pm 0.2$ (4) &\textit{0.79} & $\textbf{111.24}\pm 0.3$ (1) &\textit{0.67} & $111.07 \pm 0.1$ (3) &\textit{1.00} & $111.14 \pm 0.0$ (2) &\textit{1.00}\\ 
9  & $\textbf{111.70}\pm 0.5$ (1) &\textit{0.68} & $111.46 \pm 0.1$ (2) &\textit{0.68} & $110.19 \pm 0.7$ (3) &\textit{1.00} & $109.36 \pm 2.9$ (4) &\textit{0.98}\\ 
\bottomrule
\end{tabular*}
\vspace*{-0.4cm}
\end{table}

\begin{figure}[t]
\vspace{0.2cm}
\includegraphics[width=\textwidth]{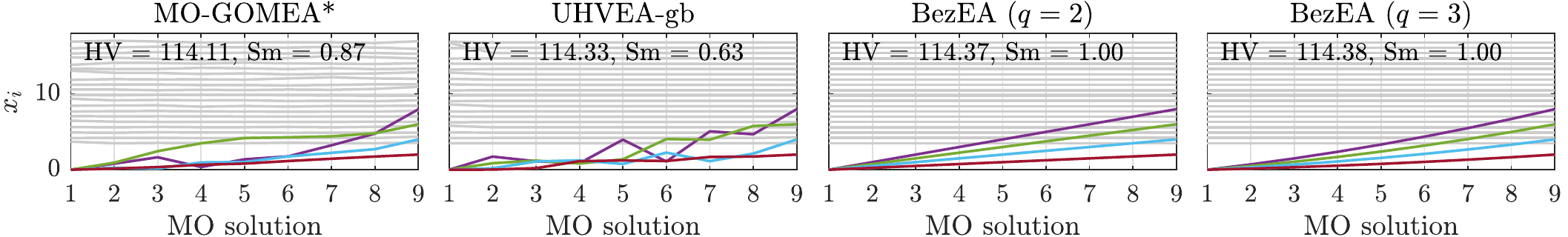}
\vspace*{-0.6cm}
\caption{Parallel coordinate plots shows of decision variables $x_i$ for WFG7. In color the $k_\text{WFG} = 4$ position-type decision variables, in grey the remaining decision variables.}
\label{fig:ppsn20_wfg_parallel}
\vspace{-0.4cm}
\end{figure}

Results are given in Table~\ref{tab:ppsn20_wfg}.
WFG1 is problematic, as none of the algorithms have an explicit mechanism to deal with its flat region. WFG2 has a disconnected Pareto front. MO-GOMEA* and UHVEA-gb both obtain solutions in multiple subsets, while BezEA obtains all solutions in a single connected subset, and spreads out well there. The linear front of WFG3 corresponds to the equidistant distribution of solutions along the B\'ezier curve, and BezEA outperforms the other methods there. Increasing $q$ generally increases performance of BezEA, except for WFG6 and WFG9. Both these problems are non-separable, and require a larger population size than the currently used $N = 200$ to be properly solved.  However, the  guideline for non-separable problems results in a population size that is too large to be of practical relevance here. In terms of smoothness, BezEA with $q = 3$ is able to obtain a smoothness close to 1, while simultaneously obtaining the best $\mbox{HV}_9$ for 4/9 problems. MO-GOMEA* obtains a mean smoothness of 0.81 while UHVEA-gb obtains the worst mean smoothness (0.68). To illustrate the obtained smoothness a parallel coordinate plot for WFG7 is given in Figure~\ref{fig:ppsn20_wfg_parallel}. This figure shows a clear pattern in decision variable values along the front (in the order $\bo$) for BezEA. This pattern is not obvious for the other two methods, while they achieve only a slightly lower hypervolume, and a lower smoothness.

\section{Discussion and outlook}
\label{sec:ppsn20_discussion}
In this work, we parameterized approximation sets as smooth B\'ezier curves in decision space, thereby explicitly enforcing a form of smoothness between decision variables of neighboring solutions when the approximation front is traversed, aimed to improve its navigability. We used an UHV-based MO problem formulation that directly allows for the optimization of parameterized approximation sets. Solving this B\'ezier problem formulation with GOMEA (BezEA), was shown to be competitive to UHV-based optimization and domination-based MOEAs, while smoothness is guaranteed. We showed that approximation sets obtained with BezEA show a more clear pattern in terms of decision variables when traversing the approximation front on a set of benchmark problems, which suggests that this approach will lead to a more intuitive and smooth approximation set navigability for real-world optimization problems.

We chose to fix the solution set size $p$ for BezEA during and after optimization, but since a parametric expression of the approximation set is available, it is straightforward to construct a large approximation set after optimization. This could be exploited to increase performance of BezEA, as it currently show computational overhead on the simple bi-sphere problem in terms of multi-objective function evaluations compared to UHVEA. In contrast to MOEAs, UHVEA and BezEA have the ability to converge to the Pareto set. When the problem is multimodal, UHVEA will spread its search over multiple modes. In that case, even an a posteriori fitting of a smooth curve through the obtained approximation set will result in low-quality solutions. BezEA on the other hand aims to obtain solutions in a single mode, thereby guaranteeing smoothness, even in a multimodal landscape. This form of regularization that is enforced upon approximation sets shows that BezEA can outperform MO-GOMEA* and UHVEA-gb on multiple problems in the WFG benchmark.


The smoothness measure introduced in this work is a measure for entire solution sets $\cS_p$, and not for individual solutions $\bx$. It can therefore not be added directly as an additional objective to the original multi-objective problem $\bff(\bx)$. We chose in this work to introduce a parameterization of approximation sets that directly enforces smoothness. Alternatively, smoothness could also be added as a second objective to the UHV-based problem formulation. This then results in the $pn$-dimensional bi-objective optimization problem, given by $h(\cS_p) = [  \text{UHV}_\bff(\cS_p) \; ; \; \text{Sm}(\cS_p) ]$. This problem can then be solved with a domination-based MOEA, or even by again formulating it as a (much) higher-dimensional UHV-based single-objective problem. Whether this approach can be efficient, even when grey-box properties such as partial evaluations are exploited, remains however future work.

The problems in this work were limited to problems involving two objectives. The presented results show that it is an interesting research avenue to extend this work to problems with more objectives. The Pareto front of non-degenerate problems with $m$ objectives is an $m-1$-dimensional manifold. Instead of a one-dimensional B\'ezier curve, the Pareto set can then be modeled by an $(m-1)$-dimensional B\'ezier simplex \cite{Kobayashi18}. For the navigation of higher-dimensional manifolds, a one-dimensional path through all obtained solutions could still be used. However, navigation would be performed might be problem specific and should be discussed with end-users. BezEA is applied to treatment planning of brachytherapy for prostate cancer, and results can be found in the  supplementary of this work. 

\textit{Source code for the algorithms in this work is made available at \cite{githubUHV}.}


\subsection*{Acknowledgments}
{\small
This work was supported by the Dutch Research Council (NWO) through Gravitation Programme Networks 024.002.003. We furthermore acknowledge financial support of the Nijbakker-Morra Foundation for a high-performance computing system.
}

\bibliographystyle{splncs04}
\bibliography{Maree_2020}

\clearpage
\section{Supplement: Solving a real-world optimization problem in brachytherapy for prostate cancer}
\label{sec:supplementary}
\label{sec:bbez20_brachy}
\beginsupplement

We demonstrate BezEA on a real-world bi-objective optimization problem that arises in the treatment of prostate cancer with brachytherapy \cite{Luong2018swarm}. Brachytherapy is a form of internal radiation therapy. In brachytherapy for prostate cancer, catheters are temporarily placed in, or close to, the prostate. Through these hollow catheters, a radioactive source can be moved, which can be stopped at predefined \textit{dwell positions}. The longer the source dwells at a certain position, the more the surrounding tissue is irradiated. The set of dwell times is called a \textit{treatment plan}. Treatment planning is the process of determining these \textit{dwell times}, such that the tumor is irradiated as much as possible, while surrounding healthy tissue is spared as much as possible. Brachytherapy treatment planning is therefore inherently a multi-objective optimization problem. 

\begin{figure}[t]
\begin{minipage}{0.49\textwidth}
\includegraphics[width=\textwidth]{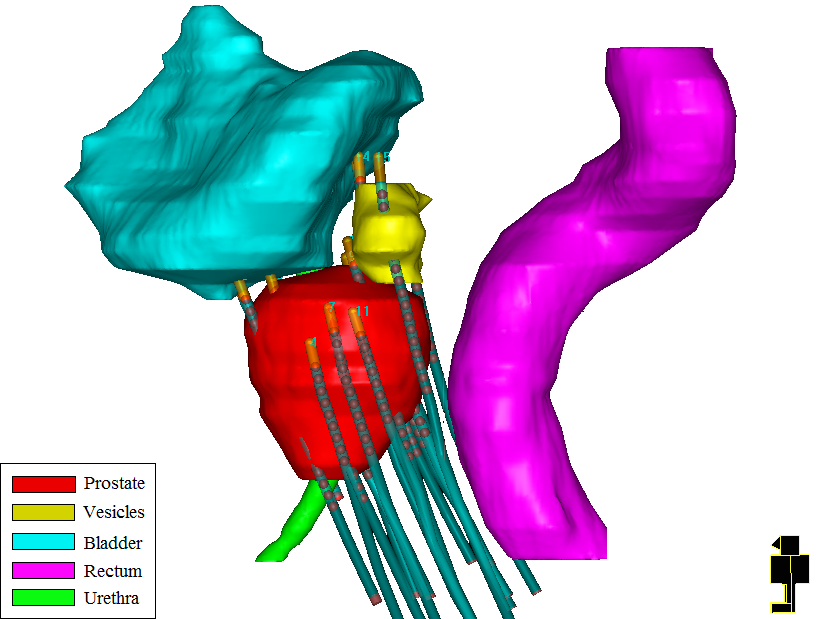}
\end{minipage}
\begin{minipage}{0.5\textwidth}
\includegraphics[width=\textwidth]{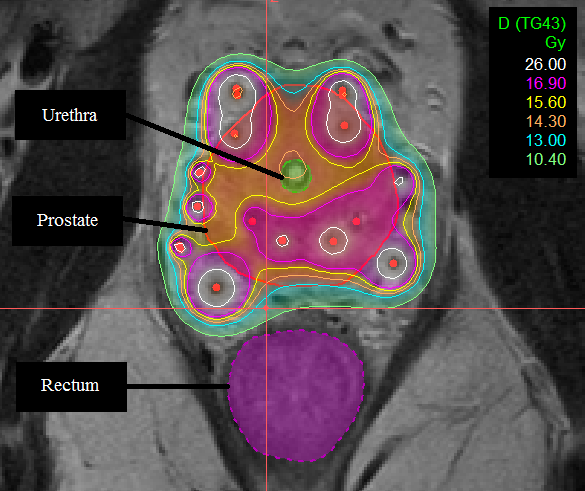}
\end{minipage}
\caption{Left: 3D reconstruction of the patient anatomy based on MR imaging, Dwell positions are indicated as spheres within each catheter. Right: MR image with delineated prostate, urethra, and rectum visible. Red dots indicate dwell positions within catheters (that are perpendicular to the screen). Thin contour lines, as well as the shaded colored heatmap, show the projected dose distribution. The aim is that the 13~Gy (light blue) contour line surrounds the prostate (in red) as much as possible, without including too much of the other tissue. Note that this MR image shows a single slice of the 3D volume.}
\label{fig:bbez20_dose_distribution}
\end{figure}

Treatment planning is performed based on magnetic resonance (MR) imaging, from which a 3D model of the patient is constructed (see Figure~\ref{fig:bbez20_dose_distribution}). For this, a radiation oncologist and radiation treatment technologist manually delineates the important structures on the MR images: the location of the catheters, and thereby the location of the dwell positions; the tumor or the \textit{target volumes}, which are the volumes that need to be treated; and other important organs and structures that need to be spared, which are referred to as the organs at risk (OARs). An example of such a 3D model of the patient is shown in Figure~\ref{fig:bbez20_dose_distribution}. Using the 3D model, the radiation \textit{dose distribution} can be simulated \cite{Rivard2004,VanderMeer2019}. The dose distribution can then be projected as a heatmap on top of the MR images, as shown in Figure~\ref{fig:bbez20_dose_distribution}, can can be used to visually determine the quality of the treatment plan. Additionally, a number of indicators of the dose distribution have been formulated, called dose-volume indices (DVIs), that were found to correlate well with different aspects of the treatment outcome, such as the five-year survival rate and severity of adverse effects for different OARs \cite{Hoskin2013}. 
The least coverage index (LCI) and least sparing index (LSI) can be used to aggregate all of these DVIs into two objective functions, resulting in a bi-objective maximization problem \cite{Luong2018swarm}. By construction, if $\text{LCI} > 0$, all coverage-related DVIs are of acceptable quality, and similarly, if $\text{LSI} > 0$, all sparing-related DVIs are of acceptable quality. The aim of bi-objective treatment planning is thus to obtain treatment plans with both LCI and LSI larger than zero (see e.g., Figure~\ref{fig:bbez20_fronts}), although this is not always achievable due to the geometric properties of the patient anatomy or the locations of the implanted catheters.


The bi-objective planning model can be solved efficiently with the multi-objective gene-pool optimal mixing evolutionary algorithm (MO-GOMEA),. It was shown that MO-GOMEA outperforms other well-known multi-objective evolutionary algorithms on this problem \cite{Luong2018swarm}. In a follow up study, it was furthermore shown that obtained treatment plans from the bi-objective planning model with MO-GOMEA were preferred in 98\% of the cases over clinically used treatment plans \cite{maree19observer}, and this approach to treatment planning results in plans of similar or better compared to other treatment planning methods and models \cite{maree20tuning}.

The main advantage of taking a bi-objective approach to treatment planning is that the resulting approximation set gives insight in the patient-specific trade-offs between coverage and sparing. Since only two objectives are used, the resulting plans can be visualized as a trade-off curve (of which we will see examples later in Figure~\ref{fig:bbez20_fronts}). It was shown that radiation oncologists appreciated the insight gained from being able to compare multiple treatment plans \cite{maree19observer}. In the selection of a desirable plan, additional patient-specific information was used, such as tumor stage, the patient's age, previous treatments or overall health. These aspects are not taken into account in the objective functions as it is not directly clear how to quantify or how to combine them. As not all information about a treatment plan is included in the objective values, a preselection of a small number of plans with desirable trade-offs in the LCI and LSI is made based on the visualized trade-off curve. The dose distribution of these preselected plans is then visually inspected, in order to select a single desirable plan. This inspection is however time consuming. If plans along the trade-off curve vary smoothly in terms of the underlying dwell times, there is an intuitive or sensible variation in properties of these plans, which will make the inspection of plans with similar trade-offs intuitive and user friendly. If dwell times of plans do not vary smoothly along the trade-off curve, it could be necessary to inspect all individual plans in order to be sure that the most desirable plan has been selected, which is time consuming and therefore infeasible to implement in clinical practice. We aim to overcome this limitation by solving the bi-objective planning model with BezEA, thereby enforcing a smoothly navigable trade-off curve. 

\subsection{Problem definition}
Let us formally define the optimization problem. Denote a treatment plan by a set of dwell times $\bt \in\bbR^n_{\geq0}$. Typically, depending on the patient and the number of implanted catheters, a few hundred dwell times need to be optimized. The computation of the dose distribution $\bd\in\bbR^{n_d}$ consists of a large matrix-vector multiplication $\bd = R\bt$, where the \textit{dose-rate matrix} $R^{n_d\times n}$ can be precomputed before optimization. We set the number of \textit{dose calculation points} $|\bd|$ to $n_d = |\bd| = 2\cdot 10^4$ during optimization \cite{Bouter2019medphysGPU}. These points are randomly sampled within the relevant structures in the 3D model, and fixed during the entire run. Using more dose calculation points during optimization increases computation time. As all calculations are based on randomly sampled dose calculation points, the obtained results naturally inhibit some uncertainty \cite{VanderMeer2019}. One could therefore recompute the dose distribution after optimization using more points (e.g., $10^6$) to remove any potential over-fitting bias from the final results before presenting plans to a clinician. 

The objective functions of the bi-objective model are formulated as $\bff^\text{brachy}(\bt) = [ \text{LCI}(R\bt) \; ; \; \text{LSI}(R\bt)]$. Note that the LCI and LSI are non-linear and non-separable functions, and both have a computational complexity of $\cO(n_d \log n_d)$. 

The aim of optimization is to obtain plans that satisfy both $\text{LCI}>0$ and $\text{LSI}>0$, but as this is not achievable for all patients, some margin is taken into account. To do so, we use a reference point, which is set to $r = (-0.04, -0.2)$.  From this, a constraint function is defined,
\begin{equation}
\label{eqn:BBez20_constraint}
C(\bt) = \max\{ -\text{LCI}(R\bt)-0.04, 0 \} +  \max\{ -\text{LSI}(R\bt)-0.2, 0 \},
\end{equation}
and constraint domination is used to handle constraint violations \cite{Deb00}. This makes sure that resulting plans always satisfy $\mbox{LCI} \geq -0.04$ and $\mbox{LSI} \geq -0.2$. The margin used here is is smaller than the margin that was used in \cite{Bouter2019medphysGPU}, as results in \cite{maree19observer} gave indication that plans with $\text{LCI} < -0.04$ were not of clinical interest. This constraint is not necessary when solving the bi-objective planning problem with UHVEA or BezEA, as simply setting the hypervolume reference point to $r$ is sufficient to guide the search.

For the computation of the objective values, it can be exploited that the required matrix-vector multiplication can be performed on a GPU \cite{Bouter2019medphysGPU}. By problem-specific tuning of MO-GOMEA, it was shown that computation time can be reduced to 30~seconds, whereas the same computations would take 2~hours on a CPU. In the same work, a form of exponential weighting of DVIs was incorporated in  the LCI and LSI, which we use here. 

\subsection{A linkage model for UHVEA and MO-GOMEA}
When only a few dwell times change, the dose distribution can be quickly updated, as only a small subset of the matrix-vector multiplication has to be performed. This is a property that MO-GOMEA can exploit to be able solve this rather high-dimensional problem with a smaller population size, and thereby in less time \cite{Luong2018swarm}. Which subsets of dwell times (i.e., decision variables) are changed simultaneously is captured in a linkage model \cite{bouter17GECCOmogomea}. To construct this linkage model, hierarchical clustering (unweighted pair group method with arithmetic mean (UPGMA)) is used to iteratively cluster dwell times together based on the distance between the corresponding dwell positions \cite{Bouter2019medphysGPU}. 

After each change of the dose distribution, even when only a few dwell times change, the LCI and LSI need to be recomputed, which gives a constant computational overhead, making it inefficient to consider small subsets of dwell times simultaneously. Therefore, the minimum number of dwell times that is changed simultaneously is set to 5 \cite{Bouter2019medphysGPU}.

To use UHVEA-gb to solve this real-world bi-objective problem, it is reformulated as a $p\cdot n$-dimensional single-objective problem, where the decision variables of $p$ MO-solutions are concatenated. For each of these MO-solutions, the same linkage model is constructed as in MO-GOMEA. All $p$ linkage models are then united into a single linkage model that is used for UHVEA-gb. Note that MO-solutions are thus updated independently (i.e., the maximum number of decision variables that are changed simultaneously is $n$). 

\subsection{B\'ezier-specific exploitable properties}
In BezEA, solution sets $\cS_{p,q} = \{\bt_1,\ldots,\bt_p\}$ of $p$ solutions are parameterized as points on a B\'ezier curve with $q$ control points. We can exploit the linearity of the dose distribution computation to reduce computation time. Intuitively, the dose distribution of each solution $\bt_i$ is a linear interpolation between the dose distributions corresponding to the control points $\bc_j$. To evaluate an entire solution set, $p$ dose distributions $\bd_i = R\bt_i$ need to be computed, from which the objective values $\bff^\text{brachy}(\bt_i) = [ \text{LCI}(\bd_i) \; ; \; \text{LSI}(\bd_i)]$ can be computed. However, since $p$ is generally larger than $q$, we use that, $$\bt_i = \bB\left(\frac{i-1}{p-1};\cC_q\right) = \sum_{j = 1}^q b_{j-1,q-1}\left(\frac{i-1}{p-1}\right)\bc_j := \sum_{j = 1}^q  b_{j,q}^{i,p} \bc_j.$$ This gives the following expression for the dose distribution computations,
$$\bd_i = R\bt_i = R \bigg( \sum_{j = 1}^q  b_{j,q}^{i,p} \bc_j \bigg)  =  \sum_{j = 1}^q  b_{j,q}^{i,p} (R \bc_j).$$
Since the product $R\bc_j$ is independent of $i$, the required number of matrix-vector multiplications reduces hereby from $p$ to $q$. Note again that the LCI and LSI are non-separable and still need to be computed $p$ times. 

Instead of directly optimizing the dwell times, $\bx = \sqrt{\bt}$ is optimized with UHVEA and BezEA, so that the search space is unbounded, which gives more freedom to BezEA to fit a curve close to $t_i = 0$, and also makes it easier for UHVEA to sample solutions close to the boundary without having to worry about boundary handling. 

\subsubsection{A linkage model for BezEA} 
For BezEA, when one decision variable of a single control point changes, the corresponding dwell time for all $p$ plans changes. It is therefore not efficient to use the same linkage model as in UHVEA, which relied on possibility to update solutions independently.
Instead, first, all corresponding decision variables of all $q$ control points are clustered together. This results in $n$ clusters each of size $q$. The same UPGMA clustering algorithm as used in MO-GOMEA is then used to construct a linkage model by iteratively merging the $n$ clusters of size $q$. In line with MO-GOMEA and UHVEA-gb, the lower bound on the cluster size is set to $5q$, such that always at least 5 dwell times of a plan are changed simultaneously. By this construction, the maximum cluster size is $qn$ (i.e., all of the decision variables).
 
\subsection{Experimental setup}
In \cite{Bouter2019medphysGPU}, a population size of $N_\text{mo} = 96$ with $K_\text{mo} = 5$ clusters was found to work well for MO-GOMEA on this problem. Using the relations $K_\text{mo} = 2p$ and $N_\text{mo} = pN$ as presented before, we deduce from this a population size of $N = 38$ for UHVEA-gb, and $N = 38q$ for BezEA. We compare BezEA with linear approximation sets ($q = 2$) to UHVEA-gb and MO-GOMEA. As it is in clinical practice only feasible time-wise to inspect a limited set of different plans, we aim the search for solution sets of size $p = 10$. For a fair and insightful comparison of MO-GOMEA with BezEA and UHVEA-gb, we again apply gHSS \cite{Guerreiro16} to reduce the obtained elitist archive (of up to 1250 plans) to an approximation set of $p$ plans for MO-GOMEA, which we denote by MO-GOMEA*.

Since all methods exploit problem-specific properties differently, we compare run time instead of MO-fevals. All methods are implemented in C++ and are run on the same CPU with a time limit of 2 hours, which corresponds to roughly 30 seconds on a GPU. BezEA and UHEA are both population-based algorithms, the entire population can be evaluated in a batch in parallel, suggesting that a similar GPU speedup can be expected as was obtained by MO-GOMEA* on a GPU.

We consider three patients here, with respectively $n = \{200, 218, 195\}$ dwell times. The third patient is known to give rise to a more difficult optimization problem, as plans with $\mbox{LCI}>0$ and $\mbox{LSI}>0$ are not achievable for this patient. To get insight in the stochastic behavior of the algorithms, we repeat all experiments 10 times and report mean and min/max performance.

\subsection{Results}

\begin{figure}[t]
\includegraphics[width=\textwidth]{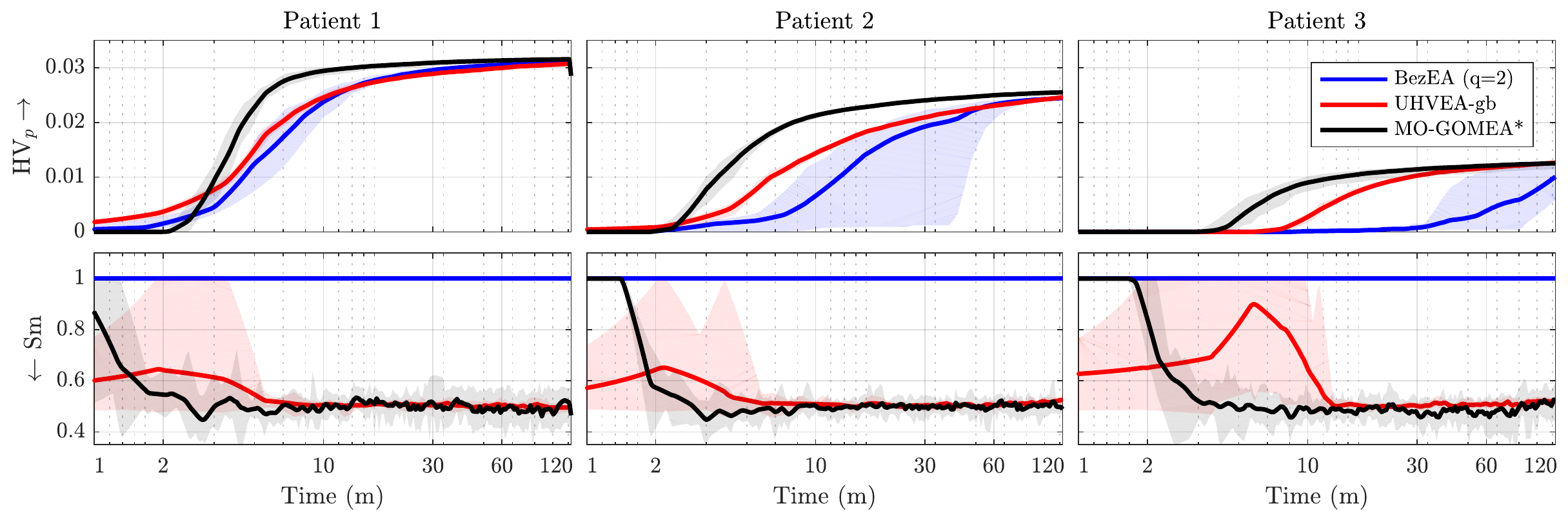}
\caption{Obtained hypervolume ($\mbox{HV}_p$) for $p = 10$ MO-solutions, and Smoothness (Sm) for three patients in terms of computation time in minutes (m). Mean values are shown over 10 runs and the shaded area shows min/max performance.}
\label{fig:bbez20_hv_smoothness}
\end{figure}

Mean $\text{HV}_p$ and smoothness results are shown in Figure~\ref{fig:bbez20_hv_smoothness}. In terms of hypervolume, all methods obtain rather similar values at the end of the run. A clear difference is observed for the difficult patient (Patient 3), where the difference between UHVEA-gb and BezEA is large initially, but ultimately, BezEA obtains similar hypervolume values. BezEA shows more variance in obtained hypervolume (as indicated by the shaded min-max performance) compared to the other two methods. This effect is largest in Patient 2 and 3, which give rise to a more difficult optimization problem (since a lower hypervolume is obtainable). As soon as a solution set is obtained that does not violate any of the constraints introduced in the B\'ezier problem formulation, the rate of convergence is however rather constant. At that point, the search is driven by hypervolume maximization, similar as in UHVEA-gb, which shows little to no variance in the obtained hypervolume values. This suggests that the constrained problem formulation of BezEA can be improved in order to better guide the search towards the feasible domain.

\begin{figure}[t]
\includegraphics[width=\textwidth]{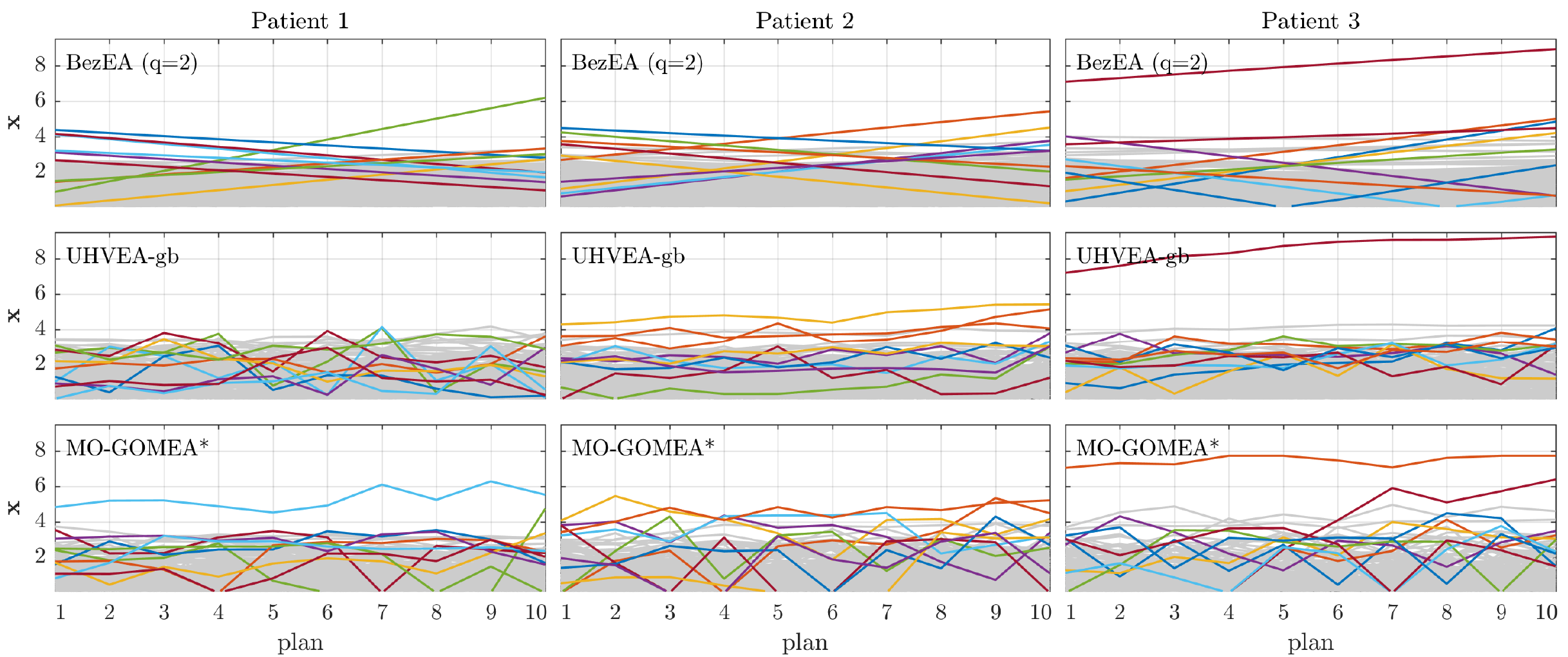}
\includegraphics[width=\textwidth]{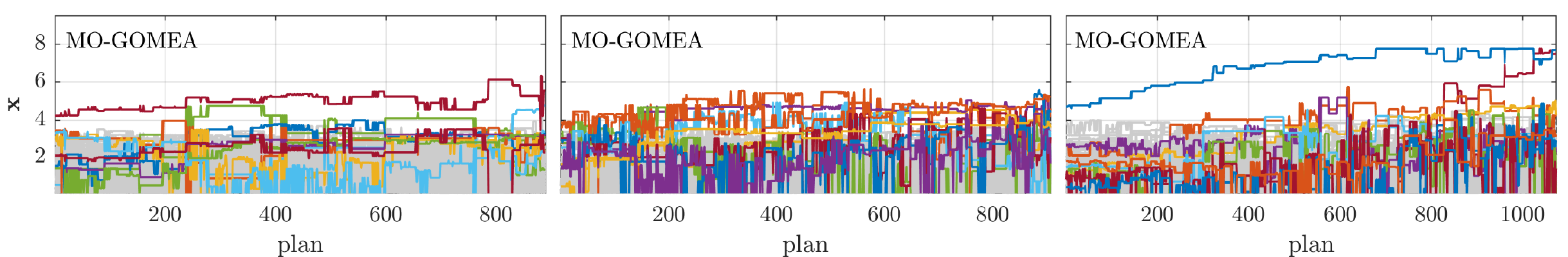}
\caption{Parallel coordinate plot of the final approximation set of a randomly selected run of each algorithm. Top three rows show the approximation set of size $p = 10$ for BezEA, UHVEA-gb and MO-GOMEA*. The bottom row shows the elitist archive of MO-GOMEA (with up to 1250 plans). On the horizontal axis, obtained plans are shown, ordered from worst to best LCI. The decision variables ($x_i =\sqrt{t_i}$, unit $\sqrt{s}$) of different plans are connect by a line. The 10 decision variables with the largest standard deviation across plans are highlighted in color, and the others are shown in grey.}
\label{fig:bbez20_parallel}
\end{figure}

BezEA with $q = 2$ has per definition a perfect smoothness of 1.0. The smoothness of both MO-GOMEA* and UHVEA-gb fluctuate around a value of 0.5, which corresponds to a zig-zag pattern in the parallel coordinate plot in Figure~\ref{fig:bbez20_parallel}. This parallel coordinate plot shows the obtained decision values. A clear difference can be observed between the approximation set obtained by BezEA and the other methods. Especially MO-GOMEA* shows a clear zig-zag pattern, indicating that plans that are next to each other on the approximation front can have very different decision values, and thereby potentially very different dose distributions.

Figure~\ref{fig:bbez20_fronts} interestingly shows that similar looking approximation fronts were obtained by the different methods, even though the corresponding decision values are very different from each other. This suggests that this real-world problem is either highly multimodal, or has (many) small plateaus.

\subsubsection{A clinical interpretation}
\begin{figure}[t]
\includegraphics[width=\textwidth]{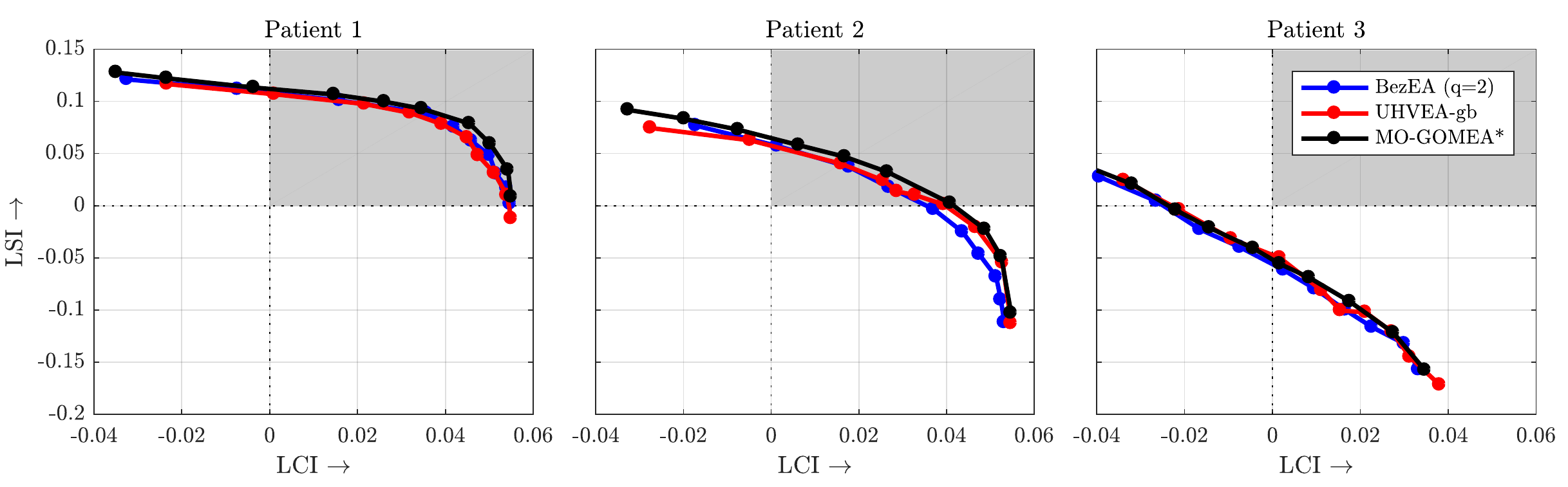}
\caption{Approximation fronts (objective space) obtained with different algorithms, of which the decision values are shown in the parallel coordinate plots in Figure~\ref{fig:bbez20_parallel}. All plans in the shaded area (with LCI$>0$ and LSI$>0$) satisfy all planning aims, and are particularly interesting from a clinical point of view. However, for some patients, it is impossible to obtain plans with these properties due to for example anatomical reasons. Patient 3 is an example of such a patient.}
\label{fig:bbez20_fronts}
\end{figure}

The objective space view (Figure~\ref{fig:bbez20_fronts}) is an important tool for decision making in clinical practice, as it gives insight in the patient-specific trade-off and maximally achievable plan quality. The obtained fronts are visually very similar, which would suggest that is does not matter which of the three fronts would be used to select a single preferred treatment plan. Furthermore, as mentioned before, especially plans with $\text{LCI} > 0$ and $\text{LSI} > 0$ are of clinical interest, if obtainable. Differences in the obtained fronts occur at the extremes of the fronts, where BezEA does not extend as far as the other methods.  The middle parts of the obtained fronts are rather similar, although BezEA obtains plans that are of slightly lower quality, but the variation in corresponding decision variables of the fronts obtained by BezEA is inherently smooth. 

To show how enforcing smoothness affects plans in the approximation set for this specific application, we analyze the LSI in more detail. The LSI is a combination of multiple planning criteria, which are used in the clinical decision making process, see \cite{maree19observer}. The volume index $V_{200\%}^{\text{prostate}}$ measures the fraction of the prostate volume that receives 200\% of the prescribed dose (13~Gy). This should not be more than 20\% of the prostate volume, which we can write as, $\Delta V_{200\%}^{\text{prostate}} = 0.2 - V_{200\%}^{\text{prostate}}$. The larger the value of $\Delta V_{200\%}^{\text{prostate}}$, the better, and the planning criteria is satisfied if $\Delta V_{200\%}^{\text{prostate}} \geq 0$. We furthermore consider the dose indices $D_{1\text{cm}^3}^{\text{bladder}}$ and $D_{1\text{cm}^3}^{\text{rectum}}$, which measure the lowest dose in the most irradiated 1~$\text{cm}^3$ of respectively the bladder and the rectum. We express the dose indices as a fraction the prescribed dose (13~Gy).§ The corresponding planning criteria are $\Delta D_{1\text{cm}^3}^{\text{bladder}} = 0.86 - D_{1\text{cm}^3}^{\text{bladder}} \geq 0$ and $\Delta D_{1\text{cm}^3}^{\text{rectum}} = 0.78 - D_{1\text{cm}^3}^{\text{rectum}} \geq 0$. The LSI combines multiple planning criteria into account by combining them in a worst-case manner, i.e.,
$\text{LSI}(R\bt) = \min\{ \Delta D_{1\text{cm}^3}^{\text{rectum}}, \Delta D_{1\text{cm}^3}^{\text{bladder}}, \ldots\}.$
By maximizing the LSI, the worst planning criteria is improved, and thus, over time, all planning criteria are improved. We refer the interested reader to \cite{Bouter2019medphysGPU} for a full description of the LSI and how exponential weighting is incorporated in this formulation. This means that the LSI takes the value of the worst planning criteria, knowing that all other planning criteria are better, but it is unknown to the optimization how much better. Plans with larger values for individual planning criteria are preferred, but by construction of the LCI, it does not directly imply that the overall fitness of these plans is also better. 

\begin{figure}[t]
\includegraphics[width=\textwidth]{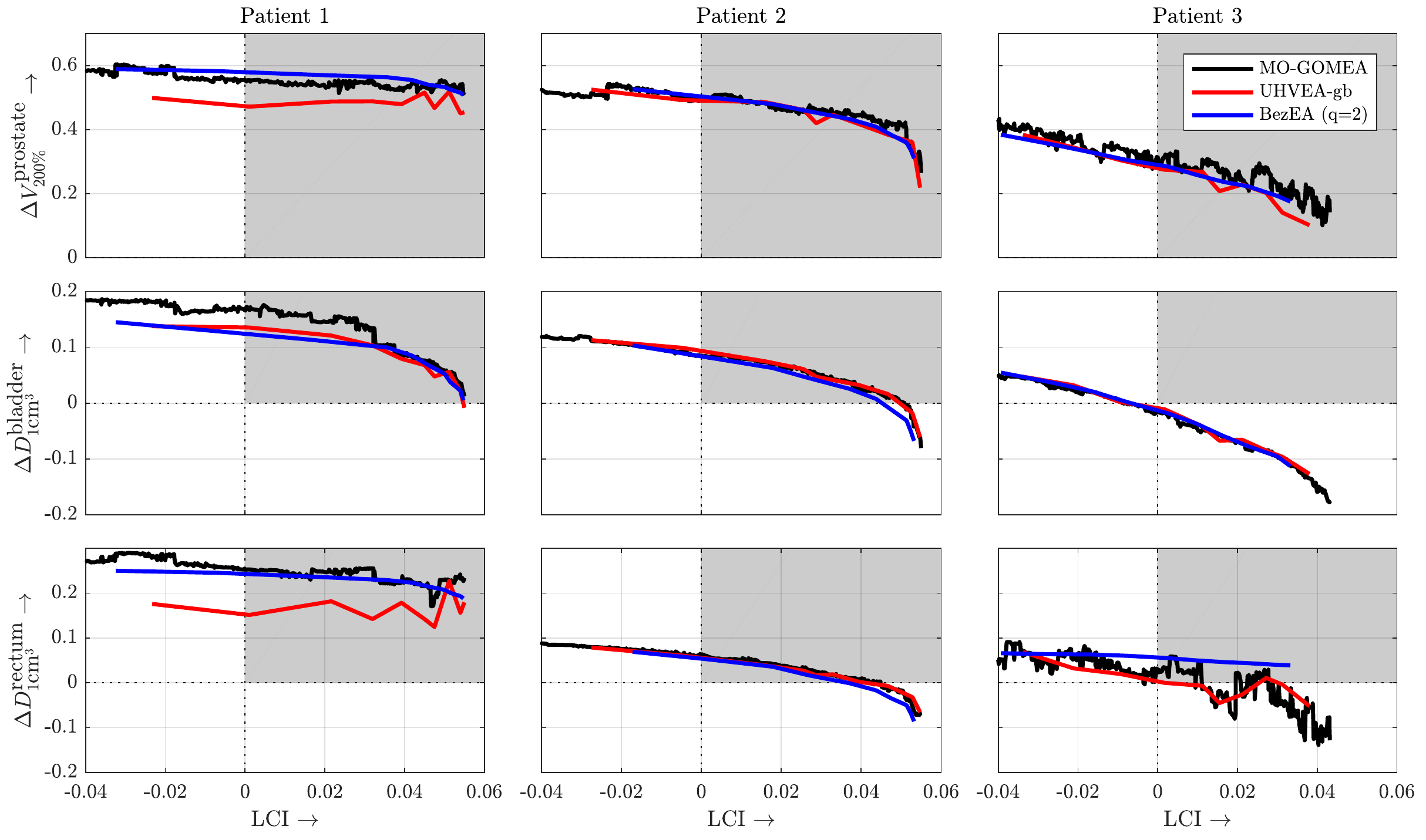}
\caption{Different DVIs, as used for clinical decision making, corresponding to the approximation sets in Figure~\ref{fig:bbez20_fronts}. For each plan, the LCI value is taken as before, but on the vertical axis, for each figure, a single DVI-based planning criterion is shown. Plans in the shaded area satisfy the coverage planning aims (i.e., LCI$>0$) and the plotted sparing criterion.}
\label{fig:bbez20_fronts_dvis}
\end{figure}

Figure~\ref{fig:bbez20_fronts_dvis} shows that values for the planning criteria obtained with BezEA fluctuate less than those obtained by with UHVEA-gb or MO-GOMEA*. For Patient 2, all algorithms obtain similar results, and for Patient 1, BezEA is slightly better than UHVEA-gb. For Patient 3, BezEA is unable to obtain plans at the right end of the front, but it obtains a slightly better rectum DVI $\Delta D_{1\text{cm}^3}^{\text{rectum}}$. Especially in terms of smoothness for Patient 3, it is clear that BezEA obtains smooth planning aim values, which suggest that the solutions on the front are more alike, which makes navigation more intuitive.

These results show that smoothly navigable approximation sets can be obtained with BezEA for the brachytherapy treatment planning problem at little to no loss in plan quality with the same computational budget as MO-GOMEA and UHVEA-gb. Differences in plan quality seem be small, suggesting that BezEA is a good alternative to MO-GOMEA while navigational smoothness is guaranteed. However, domain experts will need to be consulted to verify that plans obtained with BezEA are indeed clinically acceptable, and whether navigation of smooth approximation sets is indeed simpler and faster.

\end{document}